\documentclass{article}

\usepackage[english]{babel}

\usepackage[letterpaper,top=2cm,bottom=2cm,left=3cm,right=3cm,marginparwidth=1.75cm]{geometry}

\usepackage{amsmath}
\usepackage{amssymb}
\usepackage{amsthm}
\usepackage{tabstackengine}
\stackMath
\makeatletter
\renewcommand\TAB@delim[1]{\scriptstyle#1}
\makeatother
\setstackgap{S}{2pt}
\usepackage{xcolor}
\usepackage{graphicx}
\usepackage[colorlinks=true, hidelinks]{hyperref}
\newtheorem{theorem}{Theorem}
\newtheorem{lemma}{Lemma}
\newtheorem{corollary}{Corollary}

\begin{document}
\begin{center}
		\vskip 1cm{\LARGE \bf On the Evaluation of Apéry-Like Series Involving Multiple $t$-Harmonic Star Sums}\\
		\vskip 5mm
		Jorge Antonio González Layja\footnotemark\\
        Mexico
	\end{center}
\footnotetext{Email: \href{mailto:jorgelayja16@gmail.com}{jorgelayja16@gmail.com}}
\begin{abstract}
We evaluate, by elementary means, a new family of Apéry-like series involving multiple
$t$-harmonic star sums of even weight. Using trigonometric expansions, inverse tangent integrals, and
binomial recurrences, we obtain explicit closed-form evaluations of these series as finite alternating
sums of products of Dirichlet beta values. Several explicit examples are derived as corollaries.
\end{abstract}

\\[1ex]

\section{Introduction and Preliminaries}
In this paper, we evaluate a new family of Apéry-like series involving multiple
$t$-harmonic star sums. These sums are finite versions of the multiple $t$-star values introduced by Hoffman \cite{journal1}, which arise as odd analogues of multiple zeta-star values. Our main result provides explicit closed-form evaluations of these series as finite alternating sums of products of Dirichlet beta values. To our knowledge, such evaluations have not previously appeared in the literature.
\\[1ex]
The approach developed here is entirely elementary, relying only on trigonometric expansions, inverse tangent integrals, and binomial recurrences. In particular, we do not rely on generating-function methods, hypergeometric transformations, or complex-analytic techniques.
\\[1ex]
More precisely, for a non-negative integer $j$, we consider the Apéry-like series
$$\sum _{n=1}^{\infty }\frac{4^n}{n^2\binom{2n}{n}}t_n^{\star}(\{2\}_j),$$
where $t_n^{\star}(\{2\}_j)$ denotes the multiple
$t$-harmonic star sums of depth $j$ and weight $2j$. From this evaluation, several known and new identities follow as immediate corollaries.
\\[1ex]
Finally, we remark that the method developed in this paper adapts naturally to the corresponding multiple harmonic star sums of the same depth and weight. This allows us to recover, by different means, a conjectured identity of Gen\v{c}ev and Rucki \cite{journal2}, which was proved recently by Xu \cite{preprint1}.
\\[1ex]
We now recall the definition of the multiple
$t$-harmonic star sums of depth $j$ and weight $2j$. It is given by the nested sums
$$t_n^{\star}(\{2\}_j)=\sum _{n\ge k_1\ge \dots \ge k_j\ge 1}\prod _{i=1}^j\frac{1}{\left(2k_i-1\right)^2},$$
or equivalently by the recurrence
$$t_n^{\star}(\{2\}_j)=\sum _{k=1}^n\frac{t_k^{\star}(\{2\}_{j-1})}{\left(2k-1\right)^2},\qquad t_n^{\star}(\emptyset )=1.$$
For $j=1$, these reduce to the odd harmonic number of order $2$,
$$O_n^{(2)}=\sum _{k=1}^n\frac{1}{\left(2k-1\right)^2}.$$
More generally, the odd harmonic numbers of order $m>0$ are given by
$$O_n^{(m)}=\sum _{k=1}^n\frac{1}{\left(2k-1\right)^m}.$$
For fixed $j$, the quantities $t_n^{\star}(\{2\}_j)$ can be expressed polynomially in terms of these odd harmonic numbers. For example,
$$t_n^{\star}(\{2\}_2)=\frac{\left(O_n^{(2)}\right)^2+O_n^{(4)}}{2},\qquad t_n^{\star}(\{2\}_3)=\frac{\left(O_n^{(2)}\right)^3+3O_n^{(2)}O_n^{(4)}+2O_n^{(6)}}{6}.$$
Furthermore, we define the alternating odd harmonic numbers of order $m>0$ by
$$\overline{O}_n^{(m)}=\sum _{k=1}^n\frac{(-1)^{k-1}}{\left(2k-1\right)^m}.$$
As $n\rightarrow \infty$, $\overline{O}_n^{(m)}$ converges to the Dirichlet beta function
$$\beta (m)=\sum _{k=1}^{\infty }\frac{(-1)^{k-1}}{\left(2k-1\right)^m}.$$
The Dirichlet beta function admits explicit values at odd positive integers (see \cite[p.~807]{book2}), including
$$\beta (1)=\frac{\pi}{4},\qquad\beta (3)=\frac{\pi^3}{32},\qquad\beta (5)=\frac{5\pi^5}{1536},\qquad\beta (7)=\frac{61\pi^7}{184320}.$$
Additionally, the even value $\beta (2)=G$ is known as Catalan's constant.
\\[1ex]
We also introduce the inverse tangent integral of order $m>0$, defined for $\left|x\right|\le 1$ (see \cite[p.~190]{book4}) by
$$\operatorname{Ti}_m(x)=\sum _{k=1}^{\infty }\frac{(-1)^{k-1}x^{2k-1}}{\left(2k-1\right)^m}.$$
In particular, one has $\operatorname{Ti}_m(1)=\beta (m)$ and $\operatorname{Ti}_1(x)=\arctan (x)$.

\section{Lemmas}

\begin{lemma} \label{lma1}The following identities hold:
\begin{alignat*}{2}
(\mathrm{i})&\quad&&\text{For }m,n\in \mathbb{Z}_{\ge 0},\\
&&&\int _0^{\frac{\pi }{2}}x^{2m}\cos ((2n-1)x)\:dx=\sum _{j=0}^m\frac{\left(-1\right)^{j+n-1}}{\left(2n-1\right)^{2j+1}}\frac{(2m)!}{\left(2m-2j\right)!}\left(\frac{\pi }{2}\right)^{2m-2j}.\\
(\mathrm{ii})&\quad&&\text{For }m,k\in \mathbb{Z}_{\ge 0},\\
&&&\int _0^{\frac{\pi }{2}}x^{2m}\frac{\sin (2kx)}{\sin (x)}\:dx=2(2m)!\sum _{j=0}^m\frac{\left(-1\right)^j\overline{O}_k^{\left(2j+1\right)}}{\left(2m-2j\right)!}\left(\frac{\pi }{2}\right)^{2m-2j}.
\end{alignat*}
\end{lemma}
\begin{proof}
The result in point (i) is readily proved by repeated integration by parts. Specifically, iterating this process $2m$ times and noting the terms containing $\cos \left(\left(2n-1\right)x\right)$ vanish at both boundaries, we obtain
\begin{align*}
\int _0^{\frac{\pi }{2}}x^{2m}\cos ((2n-1)x)\:dx&=\left(\frac{\sin ((2n-1)x)}{2n-1}\frac{d^0}{dx^0}x^{2m}-\frac{\sin ((2n-1)x)}{\left(2n-1\right)^3}\frac{d^2}{dx^2}x^{2m}+\cdots \right)\bigg|_0^{\frac{\pi }{2}}\\
&=\sum _{j=0}^m\frac{(-1)^j}{\left(2n-1\right)^{2j+1}}\left(\sin ((2n-1)x)\frac{d ^{2j}}{d x^{2j}}x^{2m}\right)\bigg|_0^{\frac{\pi }{2}}.
\end{align*}
Using $\frac{d ^{2j}}{d x^{2j}}x^{2m}=\frac{(2m)!}{\left(2m-2j\right)!}x^{2m-2j}$ and applying $\sin ((2n-1)\frac{\pi }{2})=(-1)^{n-1}$ when evaluating at the endpoints, we arrive at the stated closed form.
\\[1ex]
To prove the identity in point (ii), we apply $\frac{\sin (2kx)}{\sin (x)}=2\sum _{n=1}^k\cos ((2n-1)x)$ (see \cite[p.~37]{book3}) and invoke the result in point (i). This yields
\begin{align*}
\int _0^{\frac{\pi }{2}}x^{2m}\frac{\sin (2kx)}{\sin (x)}\:dx&=2\sum _{n=1}^k\int _0^{\frac{\pi }{2}}x^{2m}\cos ((2n-1)x)\:dx.\\
&=2\sum _{n=1}^k\sum _{j=0}^m\frac{(-1)^{j+n-1}}{\left(2n-1\right)^{2j+1}}\frac{(2m)!}{\left(2m-2j\right)!}\left(\frac{\pi }{2}\right)^{2m-2j}\\
&=2(2m)!\sum _{j=0}^m\frac{(-1)^j}{\left(2m-2j\right)!}\left(\frac{\pi }{2}\right)^{2m-2j}\sum _{n=1}^k\frac{(-1)^{n-1}}{\left(2n-1\right)^{2j+1}}.
\end{align*}
Since $\sum _{n=1}^k\frac{(-1)^{n-1}}{\left(2n-1\right)^{2j+1}}=\overline{O}_k^{(2j+1)}$, this establishes the claimed identity.
\end{proof}

\begin{lemma} \label{lma2}The following identities hold:
\begin{alignat*}{2}
(\mathrm{i})&\quad&&\text{For }0<x<\frac{\pi }{2},\\
&&&\tan (x)\ln (\sin (x))=-\sum _{k=1}^{\infty }\left(\int _0^1\frac{1-t}{1+t}t^{k-1}\:dt\right)\sin (2kx).\\[1ex]
(\mathrm{ii})&\quad&&\text{For }\left|t\right|<1,j\in \mathbb{Z}_{\ge 0},\\
&&&\sum _{k=1}^{\infty }\overline{O}_k^{\left(2j+1\right)}t^{2k-1}=\frac{\operatorname{Ti}_{2j+1}(t)}{1-t^2}.\\[1ex]
(\mathrm{iii})&\quad&&\text{For }j\in \mathbb{Z}_{\ge 0},\\
&&&\int _0^1\frac{\operatorname{Ti}_{2j+1}(t)}{1+t^2}\:dt=\frac{1}{2}\sum _{k=0}^{2j}(-1)^k\beta (k+1)\beta (2j-k+1).\\[1ex]
(\mathrm{iv})&\quad&&\text{For }m\in \mathbb{Z}_{\ge 0},\\
&&&\int _0^{\frac{\pi }{2}}x^{2m}\frac{\ln (\sin (x))}{\cos (x)}\:dx\\
&&&=2(2m)!\sum _{j=0}^m\frac{(-1)^{j-1}}{\left(2m-2j\right)!}\left(\frac{\pi }{2}\right)^{2m-2j}\sum _{k=0}^{2j}(-1)^k\beta (k+1)\beta (2j-k+1).
\end{alignat*}
\end{lemma}
\begin{proof}
The identity in point $\left(\operatorname{i}\right)$ is proved in \cite[p.~243]{book1}.
\\[1ex]
To prove the result in point (ii), we note that
$$\sum _{k=1}^{\infty }\overline{O}_k^{\left(2j+1\right)}t^{2k-1}=\sum _{k=1}^{\infty }t^{2k-1}\sum _{n=1}^k\frac{(-1)^{n-1}}{\left(2n-1\right)^{2j+1}}.$$
Therefore, by interchanging the order of summation and applying the geometric series, we conclude that
$$\sum _{k=1}^{\infty }\overline{O}_k^{\left(2j+1\right)}t^{2k-1}=\sum _{n=1}^{\infty }\frac{(-1)^{n-1}}{\left(2n-1\right)^{2j+1}}\sum _{k=n}^{\infty }t^{2k-1}=\frac{1}{1-t^2}\sum _{n=1}^{\infty }\frac{(-1)^{n-1}t^{2n-1}}{\left(2n-1\right)^{2j+1}}.$$
For the identity in point (iii), integrating by parts using $\frac{d}{dt}\operatorname{Ti}_n(t)=\frac{\operatorname{Ti}_{n-1}(t)}{t}$ (see \cite[p.~190]{book4}), yields
\begin{equation*}
\int _0^1\frac{\operatorname{Ti}_{2j+1}(t)}{1+t^2}\:dt=\beta (1)\beta (2j+1)-\int _0^1\frac{\arctan (t)\operatorname{Ti}_{2j}(t)}{t}\:dt.\label{2.1}\tag{2.1}
\end{equation*}
Applying integration by parts three times to the resulting integral, we obtain
\begin{align*}
&\int _0^1\frac{\arctan (t)\operatorname{Ti}_{2j}(t)}{t}\:dt\\
&=\beta (2)\beta (2j)-\beta (3)\beta (2j-1)+\beta (4)\beta (2j-2)-\int _0^1\frac{\operatorname{Ti}_4(t)\operatorname{Ti}_{2j-3}(t)}{t}\:dt.
\end{align*}
Iterating this process $2j-1$ times, it follows that
$$\int _0^1\frac{\arctan (t)\operatorname{Ti}_{2j}(t)}{t}\:dt=\sum _{k=1}^{2j-1}(-1)^{k-1}\beta (k+1)\beta (2j-k+1)-\int _0^1\frac{\operatorname{Ti}_{2j}(t)\arctan (t)}{t}\:dt.$$
Thus,
\begin{equation*}
\begin{aligned}
\int _0^1\frac{\arctan (t)\operatorname{Ti}_{2j}(t)}{t}\:dt&=\frac{1}{2}\sum _{k=1}^{2j-1}(-1)^{k-1}\beta (k+1)\beta (2j-k+1)\\
&=\frac{1}{2}\sum _{k=1}^{2j}(-1)^{k-1}\beta (k+1)\beta (2j-k+1)+\frac{1}{2}\beta (1)\beta (2j+1).
\end{aligned}
\label{2.2}\tag{2.2}
\end{equation*}
Substituting (\ref{2.2}) in (\ref{2.1}), we arrive at
$$\int _0^1\frac{\operatorname{Ti}_{2j+1}(t)}{1+t^2}\:dt=\frac{1}{2}\beta (1)\beta (2j+1)+\frac{1}{2}\sum _{k=1}^{2j}(-1)^k\beta (k+1)\beta (2j-k+1).$$
Recognizing the trailing term as the $k=0$ term of the summation, we obtain the desired result.
\\[1ex]
To prove the result in point (iv), we consider the identity in point (ii) of Lemma \ref{lma1}. Multiplying it by $-\frac{1-t}{1+t}t^{k-1}$, integrating from $t=0$ to $t=1$, and summing from $k=1$ to $\infty$, we obtain
\begin{align*}
&\int _0^{\frac{\pi }{2}}\frac{x^{2m}}{\sin (x)}\left(-\sum _{k=1}^{\infty }\left(\int _0^1\frac{1-t}{1+t}t^{k-1}\:dt\right)\sin (2kx)\right)\:dx\\
&=2(2m)!\sum _{j=0}^m\frac{(-1)^{j-1}}{\left(2m-2j\right)!}\left(\frac{\pi }{2}\right)^{2m-2j}\int _0^1\frac{1-t}{1+t}\left(\sum _{k=1}^{\infty }\overline{O}_k^{(2j+1)}t^{k-1}\right)\:dt.
\end{align*}
Consequently, applying the substitution $t\mapsto t^2$ to the right-hand side and invoking the identities in points (i) and (ii), it follows that
$$\int _0^{\frac{\pi }{2}}x^{2m}\frac{\ln (\sin (x))}{\cos (x)}\:dx=4(2m)!\sum _{j=0}^m\frac{(-1)^{j-1}}{\left(2m-2j\right)!}\left(\frac{\pi }{2}\right)^{2m-2j}\int _0^1\frac{\operatorname{Ti}_{2j+1}(t)}{1+t^2}\:dt.$$
The stated result now follows from point (iii).
\end{proof}

\begin{lemma} \label{lma3}Let $m\in \mathbb{Z}_{\ge 0}, n\in \mathbb{Z}_{>0}$. Then the following identity holds:
$$\int _0^{\frac{\pi }{2}}x^{2m}\cos ^{2n-1}(x)\:dx=\frac{(2m)!}{2}\frac{4^n}{n\binom{2n}{n}}\sum _{j=0}^m\frac{(-1)^jt_n^{\star}(\{2\}_j)}{\left(2m-2j\right)!}\left(\frac{\pi }{2}\right)^{2m-2j}.$$
\end{lemma}
\begin{proof}
We denote the integral of interest by $I_{m,n}$. Replacing $n$ with $k$ and applying integration by parts yields
\begin{align*}
I_{m,k}&=\int _0^{\frac{\pi }{2}}x^{2m}\left(\sin (x)\right)'\cos ^{2k-2}(x)\:dx\\
&=(2k-2)\int _0^{\frac{\pi }{2}}x^{2m}\left(1-\cos ^2(x)\right)\cos ^{2k-3}(x)\:dx+\frac{2m}{2k-1}\int _0^{\frac{\pi }{2}}x^{2m-1}\left(\cos ^{2k-1}(x)\right)'\:dx\\
&=\left(2k-2\right)I_{m,k-1}-\left(2k-2\right)I_{m,k}-\frac{2m\left(2m-1\right)}{2k-1}I_{m-1,k}.
\end{align*}
This implies
$$\frac{2k-2}{2k-1}I_{m,k-1}-I_{m,k}=\frac{2m\left(2m-1\right)}{\left(2k-1\right)^2}I_{m-1,k}.$$
Multiplying both sides by $\frac{k\binom{2k}{k}}{4^k}$, using $\binom{2k}{k}=\frac{2\left(2k-1\right)}{k}\binom{2k-2}{k-1}$, and summing from $k=2$ to $n$, we obtain
$$\sum _{k=2}^n\left(\frac{\left(k-1\right)\binom{2k-2}{k-1}}{4^{k-1}}I_{m,k-1}-\frac{k\binom{2k}{k}}{4^k}I_{m,k}\right)=2m\left(2m-1\right)\sum _{k=2}^n\frac{k\binom{2k}{k}}{4^k}\frac{I_{m-1,k}}{\left(2k-1\right)^2}.$$
Hence,
$$\frac{1}{2}I_{m,1}-\frac{n\binom{2n}{n}}{4^n}I_{m,n}=2m\left(2m-1\right)\left(\sum _{k=1}^n\frac{k\binom{2k}{k}}{4^k}\frac{I_{m-1,k}}{\left(2k-1\right)^2}-\frac{1}{2}I_{m-1,1}\right).$$
Thus, applying the identity $I_{m,1}+2m\left(2m-1\right)I_{m-1,1}=\left(\frac{\pi }{2}\right)^{2m}$, which follows from two integrations by parts, we arrive at
$$I_{m,n}=\frac{1}{2}\frac{4^n}{n\binom{2n}{n}}\left(\left(\frac{\pi }{2}\right)^{2m}-2\left(2m\right)\left(2m-1\right)\sum _{k=1}^n\frac{k\binom{2k}{k}}{4^k}\frac{I_{m-1,k}}{\left(2k-1\right)^2}\right).$$
Evaluating $I_{m,n}$ for the first few values of $m$, we obtain
\begin{align*}
I_{0,n}&=\frac{1}{2}\frac{4^n}{n\binom{2n}{n}},\\
I_{1,n}&=\frac{1}{2}\frac{4^n}{n\binom{2n}{n}}\left(\left(\frac{\pi }{2}\right)^2-2\,t_n^{\star}(\{2\}_1)\right),\\
I_{2,n}&=\frac{1}{2}\frac{4^n}{n\binom{2n}{n}}\left(\left(\frac{\pi }{2}\right)^4-12\left(\frac{\pi }{2}\right)^2t_n^{\star}(\{2\}_1)+24\,t_n^{\star}(\{2\}_2)\right),\\
I_{3,n}&=\frac{1}{2}\frac{4^n}{n\binom{2n}{n}}\left(\left(\frac{\pi }{2}\right)^6-30\left(\frac{\pi }{2}\right)^4t_n^{\star}(\{2\}_1)+360\left(\frac{\pi }{2}\right)^2t_n^{\star}(\{2\}_2)-720\,t_n^{\star}(\{2\}_3)\right).
\end{align*}
The general form of the coefficients and the nested sums $t_n^{\star}(\{2\}_j)$, as suggested by these initial cases, follows from the recurrence relation by induction on $m$. This establishes the stated identity and completes the proof.
\end{proof}

\section{Main Result}
\begin{theorem} \label{thm1}Let $j\in \mathbb{Z}_{\ge 0}$. Then the following identity holds:
$$\sum _{n=1}^{\infty }\frac{4^n}{n^2\binom{2n}{n}}t_n^{\star}(\{2\}_j)=8\sum _{k=0}^{2j}(-1)^k\beta (k+1)\beta (2j-k+1).$$
\end{theorem}
\begin{proof}
By the series expansion $-\sum _{n=1}^{\infty }\frac{x^n}{n}=\ln (1-x)$ for $\left|x\right|<1$, applied with $x=\cos ^2(x)$, it follows that
$$-\frac{1}{2}\sum _{n=1}^{\infty }\frac{1}{n}\int _0^{\frac{\pi }{2}}x^{2m}\cos ^{2n-1}(x)\:dx=\int _0^{\frac{\pi }{2}}x^{2m}\frac{\ln (\sin (x))}{\cos (x)}\:dx.$$
Invoking Lemma \ref{lma3} and point (iv) of Lemma \ref{lma2}, we obtain
\begin{align*}
&\frac{(2m)!}{4}\sum _{j=0}^m\frac{(-1)^{j-1}}{\left(2m-2j\right)!}\left(\frac{\pi }{2}\right)^{2m-2j}\sum _{n=1}^{\infty }\frac{4^n}{n^2\binom{2n}{n}}t_n^{\star}(\{2\}_j)\\
&=2(2m)!\sum _{j=0}^m\frac{(-1)^{j-1}}{\left(2m-2j\right)!}\left(\frac{\pi }{2}\right)^{2m-2j}\sum _{k=0}^{2j}(-1)^k\beta (k+1)\beta (2j-k+1).
\end{align*}
Therefore, multiplying both sides of the resulting expression by $4$ and comparing the coefficients, we obtain the desired closed form.
\end{proof}

\begin{corollary}As immediate consequences of Theorem \ref{thm1}, we obtain the following evaluations:
\begin{alignat*}{2}
(\mathrm{i})&\quad&&\text{For }j=0,\\
&&&\sum _{n=1}^{\infty }\frac{4^n}{n^2\binom{2n}{n}}t_n^{\star}(\emptyset )=\sum _{n=1}^{\infty }\frac{4^n}{n^2\binom{2n}{n}}\\
&&&\hspace{6.98em}=8\beta ^2(1)=\frac{\pi ^2}{2}.\\
(\mathrm{ii})&\quad&&\text{For }j=1,\\
&&&\sum _{n=1}^{\infty }\frac{4^n}{n^2\binom{2n}{n}}t_n^{\star}(\{2\}_1)=\sum _{n=1}^{\infty }\frac{4^n}{n^2\binom{2n}{n}}O_n^{(2)}\\
&&&\hspace{8.40em}=16\beta (1)\beta (3)-8\beta ^2(2)\\
&&&\hspace{8.40em}=\frac{\pi ^4}{8}-8G^2.\\
(\mathrm{iii})&\quad&&\text{For }j=2,\\
&&&\sum _{n=1}^{\infty }\frac{4^n}{n^2\binom{2n}{n}}t_n^{\star}(\{2\}_2)=\sum _{n=1}^{\infty }\frac{4^n}{n^2\binom{2n}{n}}\frac{\left(O_n^{(2)}\right)^2+O_n^{(4)}}{2}\\
&&&\hspace{8.40em}=16\beta (1)\beta (5)-16\beta (2)\beta (4)+8\beta ^2(3)\\
&&&\hspace{8.40em}=\frac{\pi ^6}{48}-16G\beta (4).\\
(\mathrm{iv})&\quad&&\text{For }j=3,\\
&&&\sum _{n=1}^{\infty }\frac{4^n}{n^2\binom{2n}{n}}t_n^{\star}(\{2\}_3)=\sum _{n=1}^{\infty }\frac{4^n}{n^2\binom{2n}{n}}\frac{\left(O_n^{(2)}\right)^3+3O_n^{(2)}O_n^{(4)}+2O_n^{(6)}}{6}\\
&&&\hspace{8.40em}=16\beta (1)\beta (7)-16\beta (2)\beta (6)+16\beta (3)\beta (5)-8\beta ^2(4)\\
&&&\hspace{8.40em}=\frac{17\pi ^8}{5760}-16G\beta (6)-8\beta ^2(4).
\end{alignat*}
\end{corollary}

\section{Further Remarks}
The method developed in this paper adapts naturally to the following variant of the nested sums considered above:
$$\zeta_n^{\star}(\{2\}_j)=\sum _{n\ge k_1\ge \dots \ge k_j\ge 1}\prod _{i=1}^j\frac{1}{k_i^2}.$$
This represents the multiple harmonic star sums of depth $j$ and weight $2j$, which may be viewed as finite versions of the corresponding multiple zeta-star values (see \cite{journal3}).
\\[1ex]
By replacing $\cos ((2n-1)x)$ with $\cos (2nx)$ in point (i) of Lemma \ref{lma1}, the same integration-by-parts argument yields
$$\int _0^{\frac{\pi }{2}}x^{2m}\cos (2nx)\:dx=(2m)!\sum _{j=1}^m\frac{(-1)^{j+n-1}}{n^{2j}}\frac{1}{2^{2j}\left(2m-2j+1\right)!}\left(\frac{\pi }{2}\right)^{2m-2j+1}.$$
This identity is combined with the Fourier expansion $\ln (2\sin (x))=-\sum _{n=1}^{\infty }\frac{\cos (2nx)}{n}$ (see \cite[p.~45]{book3}), yielding
$$\int _0^{\frac{\pi }{2}}x^{2m}\ln (\sin (x))\:dx=(2m)!\sum _{j=0}^m\frac{(-1)^{j-1}\eta (2j+1)}{2^{2j}\left(2m-2j+1\right)!}\left(\frac{\pi }{2}\right)^{2m-2j+1},$$
where $\eta (s)$ denotes the Dirichlet eta function.
\\[1ex]
Proceeding as in the proof of Lemma \ref{lma3}, one further obtains
$$\int _0^{\frac{\pi }{2}}x^{2m}\cos ^{2n}(x)\:dx=(2m)!\frac{\binom{2n}{n}}{4^n}\sum _{j=0}^m\frac{(-1)^j\zeta_n^{\star}(\{2\}_j)}{2^{2j}\left(2m-2j+1\right)!}\left(\frac{\pi }{2}\right)^{2m-2j+1}.$$
Combining these identities as in the proof of Theorem \ref{thm1} leads to the Apéry-like series
$$\sum _{n=1}^{\infty }\frac{\binom{2n}{n}}{n4^n}\zeta_n^{\star}(\{2\}_j)=2\,\eta (2j+1).$$
This evaluation was conjectured by Gen\v{c}ev and Rucki \cite{journal2} and subsequently proved by Xu \cite[Theorem~2.2]{preprint1} using methods different from those employed here.

\end{document}